\documentclass[10pt,a4paper]{article}

\usepackage[numbers]{natbib}
\usepackage[OT1]{fontenc}
\usepackage[aop,noinfoline,amsthm,amsmath]{imsart}
\usepackage{amssymb}
\usepackage{bbm}
\usepackage{mathrsfs}
\usepackage{ae}
\usepackage{pstcol,pstricks}

\startlocaldefs
\def\@journal{Preprint}

\numberwithin{equation}{section}

\theoremstyle{plain}

\newtheorem{theorem}{Theorem}[section]
\newtheorem{lemma}[theorem]{Lemma}
\newtheorem{corollary}[theorem]{Corollary}

\theoremstyle{definition}

\newcommand{\IN}{\mathbbm{N}}
\newcommand{\IZ}{\mathbbm{Z}}
\newcommand{\IR}{\mathbbm{R}}

\def\§#1{\mathscr{#1}}
\def\£#1{\mathcal{#1}}
\newcommand{\dist}{\§L}

\def\dtv{\mathop{d_{\mathrm{TV}}}}
\def\dloc{\mathop{d_{\mathrm{loc}}}}

\newcommand{\bigo}{\mathrm{O}}
\newcommand{\lito}{\mathrm{o}}

\newcommand{\ssum}{\mathop{\textstyle\sum}}
\newcommand{\nsig}{\Sigma\kern-0.5em\raise0.2ex\hbox to 0pt{$\mid$}\kern0.5em}

\newcommand{\toinf}{\to\infty}

\newcommand{\D}{\Delta}
\newcommand{\T}{\Theta}

\newcommand{\ahalf}{{\textstyle\frac{1}{2}}}

\newcommand{\eq}{\eqref}
\newcommand{\IE}{\mathbbm{E}}
\newcommand{\IP}{\mathbbm{P}}
\newcommand{\Var}{\mathop{\mathrm{Var}}}

\newcommand{\V}{D}

\newcommand{\e}{{\mathrm{e}}}
\newcommand{\Po}{\mathop{\mathrm{Po}}}

\newcommand{\Bi}{\mathop{\mathrm{Bi}}}
\newcommand{\cBi}{\mathop{\widehat{\mathrm{Bi}}}}

\newcommand{\Be}{\mathop{\mathrm{Be}}}

\def\be#1\ee{\begin{equation*}#1\end{equation*}}
\def\ben#1\ee{\begin{equation}#1\end{equation}}

\def\bs#1\es{\begin{split}#1\end{split}}
\def\bes#1\ee{\begin{equation*}\begin{split}#1\end{split}\end{equation*}}
\def\besn#1\ee{\begin{equation}\begin{split}#1\end{split}\end{equation}}

\def\bg#1\ee{\begin{gather*}#1\end{gather*}}
\def\bgn#1\ee{\begin{gather}#1\end{gather}}

\def\bm#1\ee{\begin{multline*}#1\end{multline*}}
\def\bmn#1\ee{\begin{multline}#1\end{multline}}

\def\ba#1\ee{\begin{align*}#1\end{align*}}
\def\ban#1\ee{\begin{align}#1\end{align}}

\def\klr#1{(#1)}
\def\bklr#1{\bigl(#1\bigr)}
\def\bbklr#1{\Bigl(#1\Bigr)}
\def\bbbklr#1{\biggl(#1\biggr)}

\def\kle#1{[#1]}
\def\bkle#1{\bigl[#1\bigr]}
\def\bbkle#1{\Bigl[#1\Bigr]}
\def\bbbkle#1{\biggl[#1\biggr]}

\def\bbklel{\Bigl[}
\def\bbbklel{\biggl[}

\def\bbkler{\Bigr]}
\def\bbbkler{\biggr]}

\def\klg#1{\{#1\}}
\def\bklg#1{\bigl\{#1\bigr\}}
\def\bbklg#1{\Bigl\{#1\Bigr\}}
\def\bbbklg#1{\biggl\{#1\biggr\}}

\def\bbbklgl{\biggl\{}

\def\bbbklgr{\biggr\}}

\def\norm#1{\Vert#1\Vert}
\def\bnorm#1{\bigl\Vert#1\bigr\Vert}

\def\abs#1{\vert#1\vert}
\def\babs#1{\bigl\vert#1\bigr\vert}

\def\bbbabsl{\biggl\vert}

\def\bbbabsr{\biggr\vert}

\def\mid{\vert}
\def\bmid{\bigm\vert}

\def\floor#1{\lfloor#1\rfloor}

\def\ceil#1{\lceil#1\rceil}

\def\angle#1{\langle#1\rangle}

\def\^#1{\ifmmode {\mathaccent"705E #1} \else {\accent94 #1} \fi}
\def\~#1{\ifmmode {\mathaccent"707E #1} \else {\accent"7E #1} \fi}
\def\*#1{#1^\ast}
\def\>#1{\vec{#1}}

\def\leq{\leqslant}
\def\geq{\geqslant}

\def\atop{\@@atop}

\endlocaldefs

\begin{document}

\begin{frontmatter}
\title{Symmetric and Centered Binomial Approximation of Sums of Locally Dependent Random Variables}
\runtitle{Symmetric and Centered Binomial Approximation}

\author{%
\fnms{Adrian} \snm{R\"ollin}\ead[label=e1]{adrian.roellin@math.unizh.ch}\thanksref{t1}}
\thankstext{t1}{Partially supported by Schweizerischer Nationalfondsprojekt 20-107935/1}
\address{Universit\"at Z\"urich, Winterthurerstrasse 190, 8052 Z\"urich, Switzerland\\ \printead{e1}}
\runauthor{A. R\"ollin}

\begin{abstract} Stein's method is used to approximate sums of discrete and locally dependent random variables by a centered and symmetric Binomial distribution. Under appropriate smoothness properties of the summands, the same order of accuracy as in the Berry-Essen Theorem is achieved. The approximation of the total number of points of a point processes is also considered. The results are applied to the exceedances of the $r$-scans process and to the Mat\'ern hardcore point process type~I.
\end{abstract}


\end{frontmatter}

\section{Introduction}

The approximation of sums of dependent random variables by the standard normal distribution has been investigated in a large variety of settings. The accuracy of approximation is most often measured by the Kolmogorov and Wasserstein metrics. The use of stronger metrics typically requires that some `smoothness'-condition must be satisfied.

In this paper, under the assumption of a general local dependence structure, we study the approximation of sums of discrete random variables by a symmetric and centered Binomial distribution. This distribution serves as replacement for the normal distribution in a discrete setting. Under some general smoothness property of the summands, the same order of accuracy as in the Berry-Essen Theorem can be achieved, but now for the total variation metric. We also examine another metric, from which local limit approximations can be obtained.

In the setting of independent summands, approximation by a centered Poisson distribution has been successfully adopted by \citet{cekanavicius:01} and \citet{barbour:02}. However, for dependent summands, applications were limited to simple examples;  first attempts were made by \citet{barbour:99} and \citet{cekanavicius:01}. In contrast, the results in this paper are of general nature and allow a wide range of applications.

The proofs are based on Stein's method for distributional approximation. A main idea, introduced in \citet{ar05}, is to use interpolation functions to represent the Stein operator of a discrete distribution as the Stein operator of a continuous distribution. In the case of the Binomial, this then allows the application of standard techniques in Stein's method for normal approximation. A careful analysis of the remainder terms then shows how a suitable smoothness condition can be exploited, to obtain total variation error bounds. 

The paper is organized as follows. In the next section, we introduce the main technique in the simple case of independent summands. In section 3 these results are extended to locally dependent summands and section 4 shows their application in some examples. Section 5 contains some technical lemmas.

\subsection{Notation}

Denote by $\Bi(n,p)$ the Binomial distribution with $n$ trials of probability $p$ each. Denote by $\cBi(n,p)$ the centered Binomial distribution, i.e.\ a Binomial distribution shifted by $-n p$. Note that this distribution does not necessarily lie on the integers, but on a lattice of $\IR$ with span~$1$.

Throughout the paper, we shall be concerned with two metrics for probability distributions, the total variation metric $\dtv$ and the local limit metric $\dloc$, where, for two probability distributions $P$ and $Q$,
\ba
	\dtv\bklr{P,Q} & := \sup_{A\subset\IR}\babs{P(A)-Q(A)},\\
	\dloc\bklr{P,Q}& := \sup_{x\in\IR}\babs{P\bklr{[x,x+1)}-Q\bklr{[x,x+1)}}.
\ee
For simplicity, we will often use the notation $d_l$, where $l=1$ will stand for $\dtv$ and $l=2$ for $\dloc$.

We denote by $\norm{\cdot}$ the supremum norm if applied to functions, and the variation norm if applied to measures.
Let $\delta_x$ denote the unit mass at $x\in\IR$, and $\ast$ the convolution of measures. Define for any measure $\mu$ and any $l\in\IN := \klg{1,2,\dots}$
\be
	D^l(\mu) = \bnorm{\mu\ast(\delta_1 - \delta_0)^{\ast l}}.
\ee
Note that for measures $\mu$ and $\lambda$,
\bgn 
	\V^1\bklr{\mu} = 2\dtv\bklr{\mu,\mu\ast\delta_1},												\label{01}\\
	D^2(\mu\ast\lambda)\leq D^1(\mu)D^1(\lambda).														\label{02}
\ee

Furthermore, define $\angle{x} := x - \floor{x}$ to be the fractional part of $x\in\IR$, and $(x)_+ = x\vee 0$.

\subsection{Basic setup}

Consider a sum of the form $W=\sum_{i\in J} \xi_i$, where $W$ takes its values in a lattice of $\IR$ with span~$1$. The expectation of $W$ has no influence on the quality of the approximation, and we therefore assume without loss of generality that $\IE W = 0$; this can always be accomplished by subtracting the expectation from each individual summand. Each of the summands may now take its values on a different lattice; this, however, will result in no further complications.  

To approximate $W$ by a centered binomial distribution, we have to choose $n$ in such a way that the variance of $\cBi(n,1/2)$ is as close to the variance of $W$ as possible. As $n$ has to be integer, this is only possible up to a rounding error. However, the symmetric and centered Binomial distribution thus chosen will in general take its values on a different lattice from $W$ and the total variation distance will become $1$. To circumvent this problem, we introduce an additional parameter $t$ and approximate $W$ by a centered Binomial distribution with success probability $1/2-t$ instead ($t$ being small), to be able to match not only the variance but also the lattice.

Hence, to put the above in a rigorous form, we will make the following assumptions if not otherwise stated:

\medskip\noindent {\it Assumptions G:}\enskip
Let $J$ be a finite set and let $\{ \xi_i,i\in J\}$ be a collection of random variables with $\IE \xi_i = 0$ for all $i \in J$ and assume that there are numbers $\{a_i\in\IR;i\in J\}$ such that almost surely $\xi_i\in\IZ+a_i$. Let $W=\sum_{i\in J} \xi_i$; then $\IE W=0$ and almost surely $W\in\IZ+a$ for $a:=\sum_{i\in J} a_i$. Assume that $\sigma^2 := \Var W>1$. Define now $\delta := \angle{-4\sigma^2}$ and $t := \angle{a+2\sigma^2+\delta/2}/(4\sigma^2+\delta)$. Clearly, $4\sigma^2+\delta = \ceil{4\sigma^2}$, and by definition the distribution $\cBi\bklr{\ceil{4\sigma^2},1/2-t}$ has expectation $0$; it is also easy to check that it takes values in $\IZ+a$.

\medskip\noindent From the above definition, we see that $t$ is only of order $\bigo(\sigma^{-2})$, which is rather small in the setting that we are concerned with; Corollary \ref{co:1} shows how to obtain results without~$t$, using Lemma~\ref{lemma:a}.

\section{Sum of Independent Random Variables}

First, we examine the case of independent discrete summands. Previous work on total variation approximation has been concerned with the compound Poisson distribution (see \citet{LeCam1965} and \citet{Roos2003} and references therein), the signed compound Poisson distribution  (see \citet{Cekanavicius1997} and references therein), the Poisson distribution  (see \citet{poissonapproximation}), the centered Poisson distribution  (see \citet{MR1652516}, \citet{cekanavicius:01}, \citet{barbour:99} and \citet{barbour:02}) and some more general distributions (see \citet{Brown2001}).

We present the theorem below to demonstrate the main technique in a simple setting, noting that it also follows as a consequence of Theorem~\ref{th:2}.

\begin{theorem}\label{th:1}
Let $\{\xi_i; i\in J\}$ be independent and satisfy Assumptions G. Then, if the $\xi_i$ have finite third moments,
\be
	d_l\bklr{\dist(W),\cBi\bklr{\ceil{4\sigma^2},1/2-t}}\leq \sigma^{-2} \bbklr{\sum_{i\in J}c_{l,i}\rho_i+1.75},\qquad l=1,2,
\ee
where $\rho_i  = \sigma^3_i + \ahalf\IE\abs{\xi_i}^3$, $\sigma^2_i = \Var \xi_i$ and $c_{l,i} = \V^l\bklr{\dist(W-\xi_i)}$. 
\end{theorem}

It is clear that the above bound is useful only if the $c_{l,i}$ are small. In the case of $n$ identically distributed random variables, we need $c_{1,i}=\lito(1)$ as $n\toinf$ for asymptotic approximation in total variation, and in order to deduce a local limit theorem we must have $c_{2,i}=\lito(n^{-1/2})$. This is however always the case if $D^1(X_1)<2$ (this corresponds to the usual condition in the LLT that $X_1$ must not be concentrated on a lattice with span greater than $1$), as can be seen from \eq{55}--\eq{56}, and we then even have $c_{l,i} = \bigo(n^{-l/2})$ for $l=1,2$.

Before proving the theorem, we start with a short summary of Stein's method for Binomial approximation; for details see also \citet{stein:86} and \citet{ehm:91}. Denote by $F(M)$ the set of all real valued measurable functions on some given measure space $M$. A Stein operator $\£B:F(\IZ)\to F(\IZ)$ for the Binomial distribution $\Bi(n,p)$ is characterized by the fact that, for any integer valued random variable $W$,
\ben
	\text{$\IE (\£B g)(W) = 0$ for all bounded $g\in F(\IZ)$} \iff W\sim \Bi(n,p),						\label{03}
\ee
and a possible choice is
\ben
	(\£B g)(z) = q z g(z-1)-p(n-z)g(z),\qquad\text{for all $z\in\IZ$,}									\label{04}
\ee
where, as usual, we put  $q=1-p$. 

Let $h\in F(\IZ)$ be a bounded function. Then, the solution $g=g_h$ to the Stein equation 
\ben
	(\£B g)(z) = I[0\leq z \leq n]\bklg{h(z)-\IE h(Y)},\qquad\text{for all $z\in\IZ$,}					\label{05}
\ee
where $Y\sim\Bi(n,p)$, is also bounded. If the functions $h$ are of the form $h(z) = h_A(z) = I[z\in A]$, $A\subset\IZ$, we have the uniform bound
\ben																										\label{06}
	\norm{\Delta g_A}\leq \frac{1-p^{n+1}-q^{n+1}}{(n+1)p q},
\ee
where $\D g(z) := g(z+1) - g(z)$, and the same bound holds for $\norm{g_{\klg{b}}}$, $b\in\IZ$; see \citet{ehm:91}. Now, for all $z\in\IZ$, we can write
\be
	I[z\in A] - \IP[Y\in A] = (\£B g_A)(z) + 	I[z\notin \klg{0\dots n}]\bklr{I[z\in A]-\IP[Y\in A]},
\ee
and thus, for any integer valued random variable $V$,
\besn																										\label{07}
	\dtv\bklr{\dist(V),\Bi(n,p)} & = \sup_{A\subset\IZ} \babs{\IP[V\in A] - \IP[Y\in A]} \\
	 &\leq \sup_{A\subset\IZ}\babs{\IE(\£B g_A)(V)} + \IP\bkle{\abs{V-n/2}>n/2}.
\ee

We now construct a Stein operator for the centered Binomial distribution $\cBi(n,p)$ on the lattice $\IZ-n p$. For any function $g\in F(\IZ)$ define the function $\^g\in F(\IZ-n p)$ by $\^g(w):=g(w+n p)$ for $w\in\IZ-n p$. Then the Stein operator is defined as
\besn
	(\^{\£B}\^g)(w) & := (\£B g)(w+n p) \\ & = p(w+n p)g(w+n p)+q(w+n p)g(w-1+n p)-n p g(w+n p)\\
	& = w\bklr{p\^g(w)+q\^g(w-1)} - n p q\D\^g(w-1).					\label{08}
\ee
for all $w\in\IZ-n p$. Thus, for $W=V-n p$, an inequality corresponding to \eq{07} holds, namely
\besn
	&\dtv\bklr{\dist(W),\cBi(n,p)}\\
	&\qquad\leq \sup_{B\subset\IZ-n p}\babs{\IE(\^{\£B}\^g_B)(W)}+ \IP\bkle{\abs{W+n(p-1/2)}>n/2}. \label{09}
\ee
An equivalent inequality holds for the $\dloc$ metric, but the supremum is taken only over the sets $\klg{b}$, $b\in\IZ-n p$.

Under the assumptions of the theorem, $n=\ceil{4\sigma^2} = 4\sigma^2+\delta$ and $p=1/2-t$, and \eq{08} becomes
\ben\label{10}
	(\^{\£B}\^g)(w) = w\T\^g(w-1) - \sigma^2\D\^g(w-1) + \bklr{t^2(4\sigma^2+\delta)-wt-\delta/4}\D\^g(w-1),
\ee
where $\T\^g(w):=\ahalf\bklr{\^g(w+1)+\^g(w)}$. Since $\sigma^2>1$, the bound \eq{06} simplifies to 
\ben
	\norm{\D\^g_B} \leq \frac{1}{\sigma^2}.																\label{11}
\ee
To see this, note that $t<1/\ceil{4\sigma^2}=1/n$ and $n=\ceil{4\sigma^2}\geq5$. Then from \eq{06} we have
\be
 	\norm{\D\^g_B} \leq \frac{1}{(n+1)p q} = \frac{1}{(n+1)(1/4-t^2)}\leq \frac{4n^2}{(n+1)(n^2-4)} \leq \frac{4}{n} \leq \frac{1}{\sigma^2}.
\ee

\begin{lemma}\label{lemma:1}
Assume the conditions of Theorem \ref{th:1}. Define $\£A:F(\IZ+a)\to F(\IZ+a)$ by 
 \be
 	(\£A\^g)(w) := w\T\^g(w-1)-\sigma^2\D\^g(w-1), \qquad\text{$w\in\IZ+a$, $\^g\in F(\IZ+a)$.}
 \ee
Then,
\ben
	\babs{\IE(\£A\^g)(W)} \leq \bbklr{\norm{\D\^g}\sum_{i\in J} c_{1,i} \rho_i} \wedge \bbklr{\norm{\^g}\sum_{i\in J} c_{2,i} \rho_i}.	\label{12}
\ee
\end{lemma}

\begin{proof}
For every $w\in\IZ+a$ and $x\in[\,0,1)$ define
\ben
	f(w+x) := \T\^g(w-1)+x\D \^g(w-1) + \ahalf x^2 \D^2\^g(w-1).																\label{13}
\ee
One easily checks that $f\in C^1$ and $	f(w) = \T\^g(w-1)$ and $f'(w) = \D \^g(w-1)$, hence
\ben\label{14}
	(\£A\^g)(w) = wf(w)-\sigma^2f'(w),
\ee
for all $w\in\IZ+a$. Furthermore, $f'$ is absolutely continuous, hence $f''$ exists almost everywhere. Choose $f''$ to be the function
\ben																										 \label{15}
	f''(w+x) = \D^2\^g(w-1)
\ee
for all $w\in\IZ+a$, $0\leq x < 1$.

We can now apply the usual Taylor expansion (cf.\ \citet{reinert:98}, Theorem~2.1), but with a refined estimate of the remainder terms. Write $W_i = W-\xi_i$, $i \in J$; then
\bg
	\xi_if(W)  = \xi_i f(W_i) + \xi_i^2 f'(W_i) + \xi_i^3\int_0^1 (1-s) f''(W_i+s\xi_i)\,d s,\\
	\sigma^2_i f'(W)  = \sigma^2_i f'(W_i) + \xi_i\sigma^2_i \int_0^1 f''(W_i+s\xi_i)\,d s,
\ee
and hence, using the independence of $\xi_i$ and $W_i$ and that  $\IE \xi_i = 0$, 
\besn																										\label{16}
	\babs{\IE\bklg{\xi_if(W)-\sigma^2_if'(W)}} \leq \IE\bbbabsl \xi_i^3\int_0^1 (1-s)\IE\bkle{f''(W_i+s\xi_i)\bmid \xi_i}\,d s \\
	- \xi_i\sigma^2_i \int_0^1 \IE\bkle{f''(W_i+s\xi_i)\bmid \xi_i}\,d s \bbbabsr.
\ee
Note now that for any real valued random variable $U$ taking values on a lattice with span $1$, we obtain together with \eq{15}
\ben																										\label{17}
	\babs{\IE\bklr{f''(U+z)}} \leq  \bbklr{\norm{\D\^g}\V^1\bklr{\dist(U)}} \wedge  \bbklr{\norm{\^g}\V^2\bklr{\dist(U)}},
\ee
for all $z\in\IR$. Thus, from \eq{16} and \eq{17},
\besn
	&\babs{\IE\bklg{\xi_if(W)-\sigma^2_if'(W)}}\\
	&\leq \bbklr{\norm{\D\^g}\V^1\bklr{\dist(W_i)}\bklr{\sigma^3_i+\ahalf\IE\abs{\xi_i}^3}} \wedge
	\bbklr{\norm{\^g}\V^2\bklr{\dist(W_i)} \bklr{\sigma^3_i+\ahalf\IE\abs{\xi_i}^3}}.\label{18}
\ee
Now, using \eq{14} we have
\bes
	\babs{\IE\bklg{\£A\^g(W)}}& \leq \sum_{i\in J} \babs{\IE\bklg{\xi_if(W)-\sigma^2_if'(W)}}
\ee
and with \eq{18} the lemma is proved.
\end{proof}

\begin{proof}[Proof of Theorem \ref{th:1}]
Recall that, by Assumptions G, the distributions $\dist(W)$ and $\cBi\bklr{\ceil{4\sigma^2},1/2-t)}$ are concentrated on the same lattice. Thus, using \eq{09} and the form \eq{10} of the Stein operator, and applying the left side of the minimum in \eq{12} to the first part of \eq{10} with the bound \eq{11} gives
\besn																								\label{19}
	&\dtv\bklr{\dist(W),\cBi(4\sigma^2+\delta,1/2-t)} \\
	&\qquad \leq \frac{\sum_{i\in J} c_{1,i}\rho_i}{\sigma^2}
		+ \frac{t^2(4\sigma^2+\delta) + \sigma t+\delta/4}{\sigma^2}
		+ \IP\bkle{\abs{W}\geq 2\sigma^2-1}.
\ee
To bound the middle part of \eq{19} note that $0\leq t<(4\sigma^2+\delta)^{-1}$ and $0\leq \delta <1$. Thus, recalling that $\sigma^2 > 1$, we obtain the simple bounds
\be
	t^2(4\sigma^2+\delta) < (4\sigma^2+\delta)^{-1} \leq 1/4,
	\qquad 
	\sigma t \leq \sigma/(4\sigma^2+\delta) \leq 1/4,
	\qquad
	\delta/4 \leq 1/4.
\ee
Applying Chebyshev's inequality on the last term of \eq{19} we obtain
\be
	\IP\bkle{\abs{W}\geq 2\sigma^2-1} \leq \frac{\sigma^2}{(2\sigma^2-1)^2} \leq \frac{1}{\sigma^2}.
\ee
The $\dloc$ case is analogous, using the right side of the minimum in \eq{12} instead and the remark after \eq{06}.
\end{proof}

Note that in the next corollary we do not assume that the $\xi_i$ have expectation zero. 

\begin{corollary}\label{co:1} Let $W$ be the sum of independent and integer valued random variables $\{\xi_i,i\in J\}$ with $\sigma_i^2=\Var \xi_i$ and
\be
	v_i = \min\bklg{1/2, 1-\dtv\bklr{\dist(\xi_i),\dist(\xi_i+1)}}.
\ee
Then, if $\sigma^2>1$,
\ba
	\dtv\bklr{\dist(W),\Bi(\ceil{4\sigma^2},1/2)\ast\delta_s} &\leq 
		\frac{2\ssum\bklr{\sigma^3_i+\ahalf\IE\abs{\xi_i}^3}}{\sigma^2\klr{V - v^\ast}^{1/2}} +\frac{1+2.25\sigma^{-1}+0.25\sigma^{-2}}{\sigma}, \\
	\dloc\bklr{\dist(W),\Bi(\ceil{4\sigma^2},1/2)\ast\delta_s}& \leq
		\frac{8\ssum\bklr{\sigma^3_i+\ahalf\IE\abs{\xi_i}^3}}{\sigma^2\klr{V - 4v^\ast}_+} + \frac{3.25+0.25\sigma^{-1}}{\sigma^2},
\ee
where $s:=\ceil{\mu-\ceil{4\sigma^2}/2}$, $\mu = \IE W$, $V = \sum_{i\in J} v_i$ and $v^\ast=\max_{i\in J} v_i$.
\end{corollary}

\begin{proof} Define $W_0 = W-\mu$, and let $t$ be defined with respect to $W_0$, taking $a = -\mu$. Then, as the metrics $d_l$ are shift invariant,
\bes
	&d_l\bklr{\dist(W),\Bi(\ceil{4\sigma^2},1/2)\ast\delta_s} = d_l\bklr{\dist(W_0),\Bi(\ceil{4\sigma^2},1/2)\ast\delta_{s-\mu}}\\
	&\qquad \leq d_l\bklr{\dist(W_0),\cBi\bklr{\ceil{4\sigma^2},1/2-t}} +d_l\bklr{\Bi(\ceil{4\sigma^2},1/2-t),\Bi(\ceil{4\sigma^2},1/2)}\\&\qquad =: R^l_1 + R^l_2,
\ee
since $\Bi(\ceil{4\sigma^2},1/2-t)\ast\delta_s\ast\delta_{-\mu} = \cBi\bklr{\ceil{4\sigma^2},1/2-t}$.

Applying Lemma~\ref{lemma:a} to $R^l_2$ with the fact that $0\leq t\leq (4\sigma^2+\delta)^{-1}$ gives
\ben
	R^1_2 \leq \sigma^{-1}\bklr{1 + (2\sigma)^{-1} + (4\sigma^2)^{-1}},\qquad R^2_2 \leq \sigma^{-2}\bklr{1.5+(4\sigma)^{-1}}.			\label{20}
\ee
Define now $c_l  = \max_{i\in J}\bklg{D^l\bklr{\dist(W - \xi_i)}}$. Application of \eq{55}-\eq{56} yields
\ben
	c_1 \leq \frac{2}{\klr{V - v^\ast}^{1/2}},\qquad c_2 \leq \frac{8}{\klr{V - 4v^\ast}_+}.							\label{21}
\ee
Thus, application of Theorem \ref{th:1} to $R^l_1$ proves the corollary.
\end{proof}

\section{Locally dependent random variables}

In this section we present the main results of the paper. We exploit a finite local dependence structure as presented in \citet{MR2073183}. In the context of Stein's method for normal approximation, it has been successfully applied to a variety of problems; see for example \citet{bkr}, \citet{MR1395606} and \citet{Barbour2001}. Note that \citet{bkr} use a slightly more general dependence structure, often yielding crucial improvements when approximating sums of dissociated random variables by the normal distribution. The generalization of Theorem \ref{th:2} is straightforward, yet somewhat tedious, and we therefore use the simpler dependence structure of \citet{MR2073183}; see the Appendix for the more general version, but without proof.

Let $\{\xi_i;i\in J\}$ be a collection of random variables satisfying Assumptions G. For convenience, let $\xi_A$ denote $\{\xi_i; i\in A\}$ for every subset $A\subset J$. Assume further the following dependence structure: For every $i\in J$ there are subsets $A_i\subset B_i\subset J$ such that $\xi_i$ is independent of $\xi_{A_i^c}$, and $\xi_{A_i}$ is independent of $\xi_{B_i^c}$. Define $\eta_i = \sum_{j \in A_i} \xi_j$ and $\tau_i = \sum_{j \in B_i} \xi_j$.

\begin{theorem}\label{th:2} With $W$ as above,
\ben
	d_l\bklr{\dist(W),\cBi\bklr{\ceil{4\sigma^2},1/2-t}}\leq \sigma^{-2}\bbbklr{\sum_{i\in J}\vartheta_{l,i}+1.75},\qquad l=1,2,	\label{22}
\ee
where 
\besn																														\label{23}
	\vartheta_{l,i} ={} &
			\ahalf\IE\bklg{\abs{\xi_i}\eta_i^2 D^l\bklr{\dist(W\mid\xi_i,\eta_i)}} 
			+ \IE\bklg{\abs{\xi_i\eta_i(\tau_i-\eta_i)}D^l\bklr{\dist(W\mid\xi_i,\eta_i,\tau_i)}}\\
			& + \abs{\IE\xi_i\eta_i}\IE\bklg{\abs{\tau_i} D^l\bklr{\dist(W\mid \tau_i)}}
\ee
If further there are constants $c_{l,i}$ such that almost surely
\ben
	D^l\bklr{\dist(W\mid\xi_{B_i})}\leq c_{l,i},																				\label{24}
\ee
then
\ben
	\vartheta_{l,i} \leq c_{l,i} \bklr{\ahalf\IE\abs{\xi_i\eta_i^2} + \IE\abs{\xi_i\eta_i(\tau_i-\eta_i)} +  \abs{\IE\xi_i\eta_i}\IE\abs{\tau_i}}.
																															\label{25}
\ee
\end{theorem}
\begin{proof} Estimate \eq{25} is immediate. Following the proof of Theorem~\ref{th:1} and using Lemma~\ref{lemma:2} below, \eq{22} is proved.
\end{proof}

Note that Theorem~\ref{th:1} follows from Theorem~\ref{th:2} with the choices $A_i = B_i = \{i\}$.

\begin{lemma} 	\label{lemma:2}
Assume the conditions of Theorem~\ref{th:2}. Define $\£A:F(\IZ+a)\to F(\IZ+a)$ as in Lemma~\ref{lemma:1}. Then,
\ben
	\babs{\IE(\£A\^g)(W)} \leq \bbklr{\norm{\D\^g}\sum_{i\in J} \vartheta_{1,i}} \wedge \bbklr{\norm{\^g}\sum_{i\in J} \vartheta_{2,i}}.	\label{26}
\ee
\end{lemma}
\begin{proof} We follow the proof of Lemma \ref{lemma:1} right up to the end of the paragraph of \eq{15}. Note now that
\ben
	\sigma^2 = \sum_{i\in J} \IE\klg{\xi_i \eta_i}												\label{27}
\ee
and that, by Taylor expansion, almost surely
\besn																							\label{28}	
	\xi_i f(W) & = \xi_i f(W-\eta_i) + \xi_i\eta_i f'(W-\eta_i) + \xi_i\eta_i^2\int_0^1 f''(W-\eta_i + s\eta_i)\,d s,\\
	\xi_i\eta_i f'(W-\eta_i) &= \xi_i\eta_i f'(W-\tau_i)+\xi_i\eta_i(\tau_i-\eta_i)\int_0^1 f''(W-\eta_i + s(\tau_i-\eta_i))\,d s,\\
	\IE\klg{\xi_i\eta_i}f'(W) &= \IE\klg{\xi_i\eta_i}f'(W-\tau_i) + \IE\klg{\xi_i\eta_i}\tau_i \int_0^1 f''(W + s\tau_i)\,d s.
\ee
Now, using the facts that $\IE\xi_i = 0$, that $\xi_i$ is independent of $W-\eta_i$ and that $\eta_i$ is independent of $W-\tau_i$, we obtain from \eq{27} and \eq{28} that
\bes	
	&\IE\bklg{Wf(W)-\sigma^2f'(W)}= \sum_{i\in J}\IE\bklg{\xi_i f(W)- \IE\klg{\xi_i\eta_i}f'(W)} \\ 
	& \qquad = \sum_{i\in J} \IE\bbbklgl \xi_i \eta_i^2\int_0^1 (1-s)\IE\bklg{f''(W-\eta_i+s\eta_i)\bmid\xi_i,\eta_i}\,d s \\
	& \qquad\qquad\qquad + \xi_i\eta_i(\tau_i-\eta_i)\int_0^1 \IE\bklg{f''(W-\tau_i+s(\tau_i-\eta_i))\bmid \xi_i,\eta_i,\tau_i}\,d s \\
	& \qquad\qquad\qquad - \IE\bklg{\xi_i\eta_i}\tau_i\int_0^1\IE\bklg{f''(W-\tau_i+s\tau_i))\bmid \tau_i}\,d s  \bbbklgr.
\ee
With \eq{14} and \eq{17} the lemma follows.
\end{proof}

We now give a point process version of Theorem \ref{th:2}, exploiting mainly the same dependency structure as before.

\begin{theorem}\label{th:3}
Let\/ $\Phi$ be a simple point process on a Polish space $J$ with mean measure $\mu$. For all points $\alpha\in J$, assume that there are measurable subsets $A_\alpha\subset B_\alpha \subset J$, such that for every $\alpha\in J$
\bgn
 	\dist\bklr{\Phi_\alpha(A_\alpha^c)}=\dist\bklr{\Phi(A_\alpha^c)}, 						\label{29} \\
 	\text{$\Phi_\alpha(A_\alpha)$ and $\Phi_\alpha(B_\alpha^c)$ are independent,} 				\label{30}\\
 	\text{$\Phi(A_\alpha)$ and $\Phi(B_\alpha^c)$ are independent,}								\label{31}
\ee
where $\Phi_\alpha$ denotes the Palm process at point $\alpha$. Then, for $W=\Phi(J)-\mu(J)$ and if $\sigma^2>1$,
\besn
	& d_l\bklr{\dist(W),\cBi\bklr{\ceil{4\sigma^2},1/2-t}}\\ &\qquad \leq \sigma^{-2}\int_{\alpha\in J}\vartheta_l(\alpha) \mu(d\alpha) + 1.75\sigma^{-2},\qquad l=1,2.																	\label{32}											
\ee
where, with $\Phi':=\Phi-\mu$ and $\Phi'_\alpha:=\Phi_\alpha-\mu$,
\besn																											\label{33}
	\vartheta_l(\alpha) = {} & \babs{\IE\Phi'_\alpha(A_\alpha)}\IE\bklg{\abs{\Phi'(B_\alpha)}D^l\bklr{\Phi(B_\alpha^c)\bmid\Phi(B_\alpha)}} \\
		&+\ahalf\IE\bklg{\Phi'_\alpha(A_\alpha)^2 D^l\bklr{\Phi_\alpha(A_\alpha^c)\bmid \Phi_\alpha(A_\alpha)}} \\
		&+\ahalf\IE\bklg{\Phi'(A_\alpha)^2 D^l\bklr{\Phi(A_\alpha^c)\bmid \Phi(A_\alpha)}}\\
		&+\IE\bklg{\abs{\Phi'_\alpha(A_\alpha)\Phi'_\alpha(B_\alpha\setminus A_\alpha)}D^l\bklr{\Phi_\alpha(B_\alpha^c) \bmid \Phi_\alpha(A_\alpha), \Phi_\alpha(B_\alpha)}} \\
		&+\IE\bklg{\abs{\Phi'(A_\alpha)\Phi'(B_\alpha\setminus A_\alpha)}D^l\bklr{\Phi(B_\alpha^c) \bmid \Phi(A_\alpha), \Phi(B_\alpha)}}.
\ee
\end{theorem}
\begin{proof} Following the proof of Theorem \ref{th:1} and Lemma \ref{lemma:1}, it is clear that we only have to bound 	$\IE\bklg{Wf(W)-\sigma^2f'(W)}$ for $f$ defined as in \eq{13}.
In what follows, all integrals are taken over $\klg{\alpha\in J}$ if not otherwise stated. Note first that, because of \eq{29},
\be
	\sigma^2 = \IE\klg{\Phi(J)\Phi'(J)} = \int \mu(d\alpha)\IE\klg{\Phi_\alpha'(A_\alpha) +\Phi_\alpha'(A_\alpha^c)} = \int \mu(d\alpha)\IE \Phi_\alpha'(A_\alpha)
\ee
and hence with Taylor expansion
\bes
	\sigma^2\IE f'(W) ={}& \int\mu(d\alpha)\IE \Phi_\alpha'(A_\alpha)\IE f'(\Phi'(B_\alpha^c))\\
						& +\int\mu(d\alpha)\IE \Phi'_\alpha(A_\alpha)\IE\bbbklg{\Phi'(B_\alpha)\int_0^1f''(\Phi'(B_\alpha^c)+t\Phi'(B_\alpha))\,dt}\\
 			=: {}&  R_1 + R_2.
\ee
Now, again by Taylor, 
\bes
	\IE\klg{Wf(W)} = {}& \int \mu(d\alpha)\bkle{\IE f(\Phi_\alpha'(J))-\IE f(\Phi'(J))} \\
	 ={} &  \int\mu(d\alpha)\bkle{\IE f(\Phi'_\alpha(A_\alpha^c)) - \IE f(\Phi'(A_\alpha^c)) } \\
		& + \int\mu(d\alpha)\bbkle{\IE\bklg{\Phi'_\alpha(A_\alpha)f'(\Phi'_\alpha(A_\alpha^c))}-
			\IE\bklg{\Phi'(A_\alpha)f'(\Phi'(A_\alpha^c))}} \\
		& + \int\mu(d\alpha)\bbbkle{\IE\bbklg{\Phi'_\alpha(A_\alpha)^2\int_0^1 (1-t)f''\bklr{\Phi'_\alpha(A_\alpha^c)+t\Phi'_\alpha(A_\alpha)}\,dt}\\
	& \qquad\qquad\qquad-\IE\bbklg{\Phi'(A_\alpha)^2\int_0^1 (1-t)f''\bklr{\Phi'(A_\alpha^c)+t\Phi'(A_\alpha)}\,dt}} \\
	& =:R_3 + R_4 + R_5
\ee
and 
\bes
	R_4 = {} &\int\mu(d\alpha)\bkle{\IE\bklg{\Phi'_\alpha(A_\alpha)f'(\Phi'_\alpha(B_\alpha^c))}-\IE\bklg{\Phi'(A_\alpha)f'(\Phi'(B_\alpha^c))}}\\
		&+ \int\mu(d\alpha)\bbbklel 
			\IE\bbklg{\Phi'_\alpha(A_\alpha)\Phi'_\alpha(B_\alpha\setminus A_\alpha) \int_0^1 f''\bklr{\Phi'_\alpha(B_\alpha^c)+t\Phi'_\alpha(B_\alpha\setminus A_\alpha)}\,dt}\\
		&\qquad\qquad\qquad-\IE\bbklg{\Phi'(A_\alpha)\Phi'(B_\alpha\setminus A_\alpha) \int_0^1 f''\bklr{\Phi'(B_\alpha^c)+t\Phi'(B_\alpha\setminus A_\alpha)}\,dt}  \bbbkler \\
	 &=: R_6 + R_7.
\ee
Using \eq{29}--\eq{31}, we see that $R_3 = 0$ and $R_1 = R_6$, hence
\be
	\babs{\IE\bklg{Wf(W)-\sigma^2f'(W)}}\leq \abs{R_2}+\abs{R_5}+\abs{R_7}.
\ee
With \eq{17} we finally obtain
\ba
	\abs{R_2} &\leq \norm{\D \^g}\int\mu(d\alpha) \babs{\IE\Phi'_\alpha(A_\alpha)}\IE\bklg{\abs{\Phi'(B_\alpha)}D^1\bkle{\dist\bklr{\Phi(B_\alpha^c)\bmid\Phi(B_\alpha)}}}, \\
	\abs{R_5} &\leq \ahalf\norm{\D \^g}\int\mu(d\alpha)\bbklel \IE\bklg{\Phi'_\alpha(A_\alpha)^2 D^1\bkle{\dist\bklr{\Phi_\alpha(A_\alpha^c)\bmid \Phi_\alpha(A_\alpha)}}} \\
			  &\qquad\qquad\qquad\qquad\qquad + \IE\bklg{\Phi'(A_\alpha)^2 D^1\bkle{\dist\bklr{\Phi(A_\alpha^c)\bmid \Phi(A_\alpha)}}} \bbkler,\\
	\abs{R_7} &\leq \norm{\D \^g}\int\mu(d\alpha) \bbklel \IE\bklg{\abs{\Phi'_\alpha(A_\alpha)\Phi'_\alpha(B_\alpha\setminus A_\alpha)}D^1\bkle{\dist\bklr{\Phi_\alpha(B_\alpha^c) \bmid \Phi_\alpha(A_\alpha), \Phi_\alpha(B_\alpha)}}} \\
			&\qquad\qquad\qquad\qquad +\IE\bklg{\abs{\Phi'(A_\alpha)\Phi'(B_\alpha\setminus A_\alpha)}D^1\bkle{\dist\bklr{\Phi(B_\alpha^c) \bmid \Phi(A_\alpha), \Phi(B_\alpha)}}}\bbkler.
\ee

To obtain $\vartheta_2$, just replace $\norm{\D g}$ by $\norm{g}$ and $D^1$ by $D^2$ in the above bounds .
\end{proof}

\begin{corollary}\label{co:2} 
Let $\Phi$ be a simple point process satisfying \eq{29}--\eq{31}. If there is further a function $c_l(\alpha)$, such that for $\mu$-almost every $\alpha\in J$ almost surely
\ben
	D^l\bkle{\dist\bklr{\Phi(J)\bmid \Phi|_{B_\alpha}}},D^l\bkle{\dist\bklr{\Phi_\alpha(J)\bmid \Phi_\alpha|_{B_\alpha}}} \leq c_l(\alpha), \qquad l=1,2,							\label{35}
\ee
then \eq{33} satisfies
\besn																									\label{36}
	\vartheta_l(\alpha) \leq {} & c_l(\alpha)\bbklel\babs{\IE\Phi'_\alpha(A_\alpha)}\IE\abs{\Phi'(B_\alpha)} + \ahalf\IE\Phi'_\alpha(A_\alpha)^2 
		+\ahalf\IE\Phi'(A_\alpha)^2\\
		&\qquad\qquad+\IE\abs{\Phi'_\alpha(A_\alpha)\Phi'_\alpha(B_\alpha\setminus A_\alpha)} +\IE\abs{\Phi'(A_\alpha)\Phi'(B_\alpha\setminus A_\alpha)}\bbkler\\
	\leq{}&c_l(\alpha)\bbklel 1.5\IE\bklg{\Phi_\alpha(A_\alpha)\Phi_\alpha(B_\alpha)}+1.5\IE\bklg{\Phi(A_\alpha)\Phi(B_\alpha)} \\
		&\qquad\qquad+ 6\mu(A_\alpha)\mu(B_\alpha) + 4\mu(B_\alpha)\IE\Phi_\alpha(B_\alpha)\bbkler.
\ee
\end{corollary}

\section{Applications}
In what follows, we calculate only rough bounds, leaving much scope for improvement. In particular, we replace the moments in the estimates by almost sure bounds.

\subsection{Exceedances of the $r$-scans process}

We follow the notation of \citet{MR1161058}. Let $X_1,X_2,\dots,X_{n+r-1}$ be independent and identically distributed random variables with distribution function $F$. Define the $r$-scan process $R_i = \sum_{k=0}^{r-1} X_{i+k}$, $i=1,2,\dots,n$ and further $W^-_i = I[R_i\leq a]$ for $a\in\IR$. We are interested in the number $N^-=\sum_{i=1}^n W^-_i$, that is the number of $R_i$ not exceeding~$a$. With $p=\IE W^-_i = \IP[R_1\leq a]$, we have $\IE N^- = n p$ and 
\ben
	\sigma^2 = \Var W = n p\bbklr{1-p+2\sum_{d=1}^{r-1}(1-d/n)\psi(d)},									\label{39}
\ee
where $\psi(d) = \IP[R_{d+1}\leq a\mid R_1<a]-p \geq 0$.

Poisson approximations for the $r$-scan process have been extensively studied by \citet{MR1161058}. Normal approximation has been considered by \citet{MR1395606}; in particular they show, that, for fixed $r$ and $a$, $N^-$ converges in the Kolmogorov metric to the normal distribution with rate $\bigo(n^{-1/2})$. In the next theorem we achieve the same rate in total variation, and also a rate for the corresponding local limit approximation.

\begin{theorem} Assume that $F$ is continuous, $F(0)=0$, and $0\leq F(x)<F(y)$ for all $x<y$, and let $a>0$ be fixed. Then, for all $n$ such that $\sigma^2>1$,
\be
	d_l\bklr{\dist(N^--n p),\cBi\bklr{\ceil{4\sigma^2},1/2-t}} \leq C_l n^{-l/2}, \qquad l=1,2,
\ee
where the constants $C_1$ and $C_2$ are independent of $n$ and can be extracted from the proof.
\end{theorem}
\begin{proof}
We apply Theorem \ref{th:2} for $W = \sum_{i=1}^n \xi_i = \sum_{i=1}^n(W_i^- - p)$. We can set
\ba
	A_i & = \{i-r+1,\dots,i+r-1\}\cap \{1,\dots,n\},\\
	B_i & = \{i-2r+2,\dots,i+2r-2\}\cap \{1,\dots,n\}.
\ee
Then, as $\abs{A_i}\leq2r-1$, $\abs{B_i}\leq4r-3$ and $\abs{B_i\setminus A_i}\leq2r-2$, the following rough bounds are obvious:
\bg
	\IE\abs{\xi_i\eta_i^2}\leq (2r-1)^2,\quad \IE\abs{\xi_i\eta_i(\tau_i-\eta_i)}\leq (2r-1)(2r-2),\\
  	\abs{\IE\xi_i\eta_i}\IE\abs{\tau_i} \leq (2r-1)(4r-3),
\ee
thus
\ben
	\vartheta_{l,i} \leq c_{l,i}\bklr{16r^2-20r+6}							\label{40}
\ee
Consider now the block $B_1 = \sum_{i=1}^{3r-2} W^-_i$, and assume that the values $\partial B_1 = (X_1,\dots,X_{r-1})$ and $\partial B_2 = (X_{3r-1},\dots,X_{4r-2})$ are given. Define the events
\bg
	\£A :=  \bklg{a/r < X_r,\dots,X_{2r-2},X_{2r+1},\dots,X_{3r-2}\leq a(r+1)/r^2,\, 0<X_{2r}\leq a/(2r^2) }\\
	\£A_0 := \bklg{a/r< X_{2r-1} \leq a(r+1)/r^2},\quad \£A_1 := \bklg{0<X_{2r-1}\leq a/(2r^2)}.
\ee
Due to the conditions on $F$ and independence it is clear that $p_j := \IP[\£A\cap\£A_j]>0$ for $j=0,1$. Note now that
\be
	R_r = \sum_{i=r}^{2r-1} X_i > a\quad\text{on $\£A\cap\£A_0$}, \qquad R_r <a \quad\text{on $\£A\cap\£A_1$}.
\ee
Note further that $	R_s < a $ for all $s=r+1,\dots,2r-1$ on $\£A\cap(\£A_0\cup\£A_1)$. Hence
\be	
	\sum_{i=r}^{2r-1} W^-_i = r-1\quad\text{on $\£A\cap\£A_0$},\qquad\sum_{i=r}^{2r-1} W^-_i = r\quad\text{on $\£A\cap\£A_1$}.
\ee
It easy to see now by a coupling argument that
\be
	\ahalf D^1\bklr{\dist(B_1)}\leq 1-(p_0\wedge p_1) < 1.
\ee
Noting that by sequentially stringing together blocks like $B_1$, we can have $m:=\floor{n/(3r-2)}$ such blocks, which are independent given all the borders $\partial B_i$. Furthermore, for every $i$, the $R_j$ in $B_i$ depend on the $X_k$ of at most two such blocks. Therefore, defining $Z=(\partial B_1,\dots,\partial B_m)$ and using \eq{57} and \eq{58},
\ba
	D^1\bklr{\dist(W|\partial B_i, i=1,\dots,m)} &\leq \frac{2}{\bklr{\min\{1/2,p_0,p_1\}(m-2)}^{1/2}}=: c_{1,i}, \\
	D^2\bklr{\dist(W|\partial B_i, i=1,\dots,m)} &\leq \frac{8}{\min\{1/2,p_0,p_1\}(m-4)_+}=: c_{2,i}.
\ee
Clearly, $c_{l,i} = \bigo(n^{-l/2})$. Hence, putting this, \eq{39} and \eq{40} into \eq{22}, the theorem follows.
\end{proof}

\subsection{Mat\'ern hard-core process type I}

We approximate the total number of points of the Mat\'ern hard-core process type I introduced by \citet{MR0169346}. We use rectangular instead of the usual circular neighborhoods. Let $\Phi$ be the process on the $d$-dimensional cube $J = [0,1)^d\subset\IR^d$ defined as
\be
 	\Phi(B) = \sum_{i=1}^\tau I[X_i\in B] I[\text{$X_j\notin K_r(X_i)$ for all $j=1,\dots,\tau$, $j\neq i$}],
\ee
where $\tau\sim\Po(\lambda)$ and $\{X_i$; $i\in\IN\}$ is a sequence of independent and uniformly distributed random variables on $J$ and where, for $x=(x_1,\dots,x_d)\in J$  and $r>0$, $K_r(x)$ denotes the $d$-dimensional closed cube with center $x$ and side length $r$. To avoid edge effects, we treat $J$ as a $d$-dimensional torus, thus identifying any point outside $J$ by the point in $J$ which results in coordinate-wise shifting by~$1$. The process $\Phi$ is thus a thinned Poisson point process with rate $\lambda$ having all points deleted which contain another point in their $K_r$ neighborhood. For the mean measure $\mu$ of $\Phi$ we obtain
\ben
	\frac{d\mu(x)}{d x} = \lambda\e^{-c}.																	\label{41}	
\ee

We are now interested in the distribution of $\Phi(B)$ when $r$ is small an $\lambda$ large.

\begin{theorem}
Put $W:=\Phi(J)-\mu(J)$ and let $a>0$ be a fixed real number. Then, for every $\lambda$ and $r$ such that $\lambda r^d = a$ and $\sigma^2:=\Var W > 1$,
\be
 	d_l\bklr{\dist(\Phi(J)-\mu(J)),\cBi\bklr{\ceil{4\sigma^2},1/2-t}} \leq C_l \lambda^{-l/2}, \qquad l=1,2,
\ee
for constants $C_1$ and $C_2$ which are independent of $\lambda$ and can be extracted from the proof.
\end{theorem}
 
\begin{proof} We apply Corollary \ref{co:2}. We can take $A_x = K_{2r}(x)$ and $B_{x} = K_{4r}(x)$ and check that the conditions \eq{29}--\eq{31} are fulfilled. Some calculations show that the reduced second factorial moment measure $M$ satisfies
\be
	\frac{dM(x)}{d x} = \begin{cases}
							0 													& \text{if $x\in K_r(0)$,}\\
							\lambda^2\e^{-\lambda\abs{K_r(0)\cup K_r(x)}}	& \text{if $x\in K_{2r}(0) \setminus K_r(0)$,}\\
							\lambda^2\e^{-2a}	& \text{if $x\notin K_{2r}(0)$,}\\
	                      \end{cases}
\ee
compare with \citet[pp.~367, 373]{Daley1988}. Thus, $M(J)\geq\lambda^2\e^{-2a}(1-r^d)$ and
\ben
	\sigma^2 = \lambda\e^{-a} + M(J) -\mu(J)^2 \geq \lambda \e^{-a}\klr{1-a\e^{-a}}.						\label{42}
\ee
Since we can have at most $7^d$ points of $\Phi$ in $B_x$, we obtain from \eq{36} the rough estimate
\ben
	\vartheta_l(x) \leq 26\cdot 7^d c_l(x),		\label{43}
\ee
where $c_l(\cdot)$ is as in \eq{35}. To estimate $c_l(x)$ write $K_{r} = K_{r}(0)$. We have
\ba
	\IP\bkle{\Phi(K_{r}) = 0 \bmid \Phi_{K_{7r}\setminus K_{5r}}} &
		\geq \Po\bklr{\lambda\abs{K_{r}}}\klg{0} = \e^{-a} =: p_0, \\
	\IP\bkle{\Phi(K_{r}) = 1 \bmid \Phi_{K_{7r}\setminus K_{5r}}} &
		\geq \Po\bklr{\lambda\abs{K_{3r}\setminus K_{r}}}\klg{0}\cdot\Po\bklr{ \lambda\abs{K_r}}\klg{1} = a\e^{-3^d a}=: p_1.
\ee
Hence, by a coupling argument,
\ben
	\ahalf D^1\bklr{\dist(\Phi(K_{6r})|\Phi_{K_{7r}\setminus K_{5r}})}\leq 1-(p_0\wedge p_1) < 1.			\label{44}
\ee
Let now $x$ be arbitrary. Divide the space $J$ into boxes of side length $6r$, centered around $x$ (see Figure \ref{fig01}). With $m:=\floor{1/(6r)}$, we can have $m^d$ such boxes plus a remainder. Denote this remainder by $J^R$ and denote by $x_l$, $l\in\{1,\dots,m\}^d=:\£M$ the centers of the boxes where $x_{1,\dots,1}=x$. Note now that, given $\Phi$ on all the borders $K_{6r}(x_l)\setminus K_{5r}(x_l)$, $l\in \£M$ (grey area in Figure \ref{fig01}), the random variables $\Phi(K_{6r}(x_l))$, $l\in \£M$, are independent and satisfy inequality \eq{44}. Furthermore, $\Phi|_{J\setminus K_{6r}(x)}$ is independent of $\Phi|_{B_x}$, and therefore, defining $Z=\bklr{(\Phi|_{K_{6r}(x_l)\setminus K_{5r}(x_l)})_{l\in\£M}, \Phi|_{J^R}}$ and using \eq{57} and \eq{58}, we obtain
\ban
	D^1\bklr{\dist\bklr{\Phi(J)\bmid\Phi|_{B_x}}} &\leq \frac{2}{\bklr{\min\{1/2,p_0,p_1\}(m^d-2)_+}^{1/2}}=: c_1(x), \label{45}\\
	D^2\bklr{\dist\bklr{\Phi(J)\bmid\Phi|_{B_x}}} &\leq \frac{8}{\min\{1/2,p_0,p_1\}(m^d-3)_+}=: c_2(x). \label{46}
\ee
Noting that almost surely $\dist\bkle{\Phi_x\bklr{J\setminus K_{6r}(x)}\bmid\Phi_x|_{B_x}} = \dist\bkle{\Phi\bklr{J\setminus K_{6r}(x)}}$, we see that \eq{45} and \eq{46} hold also for $\Phi_x$, thus $c_l(x)$ satisfies \eq{35}. Now, recalling that $a=\lambda r^d$ is constant, we have $c_l(x) = \bigo(\lambda^{-l/2})$. Hence, putting this and \eq{41}--\eq{43} into \eq{32}, the theorem follows.
\begin{figure}

\psset{xunit=5mm,yunit=5mm}
\definecolor{nbr}{rgb}{0.93,0.93,0.93}
\begin{pspicture}(-7,-5.4)(7,5)
\psframe[linestyle=none,fillstyle=solid,fillcolor=lightgray](-6.98,-4.49)(6.98,4.99)
\psframe[linestyle=none,fillstyle=solid](-2.5,-2.5)(2.5,2.5)
\psframe[linestyle=none,fillstyle=solid](-2.5,3.5)(2.5,5)
\psframe[linestyle=none,fillstyle=solid](-2.5,-4.5)(2.5,-3.5)
\psframe[linestyle=none,fillstyle=solid](-7,-2.5)(-3.5,2.5)
\psframe[linestyle=none,fillstyle=solid](-7,3.5)(-3.5,5)
\psframe[linestyle=none,fillstyle=solid](-7,-4.5)(-3.5,-3.5)
\psframe[linestyle=none,fillstyle=solid](3.5,-2.5)(7,2.5)
\psframe[linestyle=none,fillstyle=solid](3.5,3.5)(7,5)
\psframe[linestyle=none,fillstyle=solid](3.5,-4.5)(7,-3.5)

\psframe[linestyle=none,fillstyle=solid,fillcolor=nbr](-1.5,-1.5)(1.5,1.5)

\psline(-7,+3)(+7,+3)
\psline(-7,-3)(+7,-3)

\psline(-3,-4.5)(-3,+5)
\psline(+3,-4.5)(+3,+5)

\rput(-6.853, 0.335){\footnotesize $+$}
\rput(-6.381,-3.521){\footnotesize $+$}
\rput(-6.359, 1.290){\footnotesize $+$}
\rput(-6.041,-2.375){\footnotesize $+$}
\rput(-5.278, 4.976){\footnotesize $+$}
\rput(-4.776,-3.737){\footnotesize $+$}
\rput(-4.322, 0.562){\footnotesize $+$}
\rput(-4.269, 2.013){\footnotesize $+$}
\rput(-3.589,-2.915){\footnotesize $+$}
\rput(-3.329,-0.709){\footnotesize $+$}
\rput(-3.199, 3.650){\footnotesize $+$}
\rput(-2.396, 4.316){\footnotesize $+$}
\rput(-2.010, 1.576){\footnotesize $+$}
\rput(-1.791,-3.704){\footnotesize $+$}
\rput(-1.777, 3.111){\footnotesize $+$}
\rput(-1.245, 0.135){\footnotesize $+$}
\rput(-0.801, 2.267){\footnotesize $+$}
\rput(-0.470, 4.312){\footnotesize $+$}
\rput( 0.002,-1.712){\footnotesize $+$}
\rput( 0.114,-4.091){\footnotesize $+$}
\rput( 0.356,-3.303){\footnotesize $+$}
\rput( 0.827,-0.364){\footnotesize $+$}
\rput( 1.049, 1.044){\footnotesize $+$}
\rput( 1.681,-3.774){\footnotesize $+$}
\rput( 1.740,-2.184){\footnotesize $+$}
\rput( 3.663,-3.388){\footnotesize $+$}
\rput( 3.761,-1.408){\footnotesize $+$}
\rput( 3.850, 4.913){\footnotesize $+$}
\rput( 3.863, 0.090){\footnotesize $+$}
\rput( 4.243, 1.624){\footnotesize $+$}
\rput( 4.449,-4.115){\footnotesize $+$}
\rput( 5.000,-0.863){\footnotesize $+$}
\rput( 5.055, 2.881){\footnotesize $+$}
\rput( 5.134, 0.368){\footnotesize $+$}
\rput( 5.732, 2.235){\footnotesize $+$}
\rput( 5.902, 3.967){\footnotesize $+$}
\rput( 6.408, 1.773){\footnotesize $+$}
\rput( 6.615,-3.504){\footnotesize $+$}

\psframe[linewidth=0.01](-0.5,-0.5)(0.5,0.5)

\uput{2.5pt}[65](0,0){$\scriptstyle x$}
\uput{1pt}[45](-1.5,-1.5){$\scriptstyle K_{3r}(x)$}
\rput(0,-5.5){$\underbrace{\hbox{\kern30mm}}_{\textstyle 6r}$}
\rput(0,4.5){$\overbrace{\hbox{\kern24mm}}^{\textstyle 5r}$}
\rput(1.3,2.35){$K_{6r}(x)$}
\psdot[linewidth=0.01](0,0)
\psdot[linewidth=0.01](-6,0)
\psdot[linewidth=0.01](6,0)
\end{pspicture}
\caption{\label{fig01}Matérn hard-core process type I: Given that the process $\Phi$ is known on the borders $\cup_{l\in\£M} K_{6r}(x_l)\setminus K_{5r}(x_l)$ (grey area), the boxes $\Phi|_{K_{6r}(x_l)}$, $l\in\£M$, are independent.}
\end{figure}
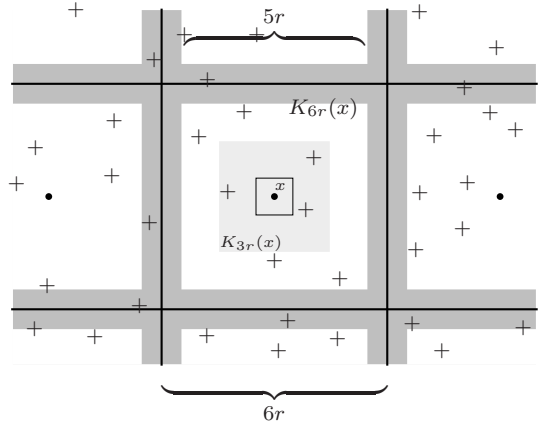
\end{proof}

\newpage


\section{Appendix}

\subsection{Properties of the solution to the Stein equation}
\begin{lemma}\label{lemma:b} For any indicator function $h(z)=I[z\in A]$, $z\in\IZ$, $A\subset\IZ$, the solution $g=g_h$ to the Stein equation \eq{05} satisfies
\ben
	\norm{g} \leq 1\wedge (n p q)^{-1/2}.																			\label{47}
\ee
\end{lemma}
\begin{proof} We apply the generator method introduced by \citet{barbour:88}. For any function $f:\bklg{0,\dots,n}\to\IR$, define
\ben
	(\£Af)(z) =  \bklr{\£B(-\Delta f)}(z) = qzf(z-1) - \bklr{qz+p(n-z)}f(z) + p(n-z)f(z+1),						\label{48}
\ee
which is the infinitesimal generator of a pure jump Markov process. A solution $g$ to \eq{05} is now given through
\be
	\psi(z) = -\int_0^\infty \IE\bklg{h\bklr{Y_z(t)}-h(Y)}\,dt,\qquad\text{for $z\in\klg{0,\dots,n},$}
\ee
and $g(z) = -\Delta \psi(z)$ for $z\in\klg{0,\dots,n-1}$ and $g(z)=0$ else, where $Y_z$ is a Markov process with generator $\£A$ starting at point $z$, and $Y$ is a random variable having the stationary distribution $\Bi(n,p)$. Now, we have for $z\in\klg{0,\dots,n-1}$,
\ben
	\Delta\psi(z) = \int_0^\infty \IE\bklg{h\bklr{Y_{z}(t)}-h\bklr{Y_{z+1}(t)}}\,dt.						\label{49}
\ee
We now fix $z$ and construct a coupling of $Y_{z}$ and $Y_{z+1}$ to bound \eq{49}. Let thereto $X^{(i)}_k(t)$, $k\in\klg{1,\dots,n}$, $i\in\klg{0,1}$, be independent Markov processes with state space $\klg{0,1}$, starting in point $i$ and having jump rate $p$ if the process is in $0$ and $q$ otherwise. It is easy to see by the Kolmogorov differential equations that 
\ben
	X^{(1)}_k(t) \sim \Be(p+q\e^{-t}),\qquad X^{(0)}_k(t) \sim \Be(p-p\e^{-t})										\label{50}
\ee
where $\Be(p)$ denotes the Bernoulli distribution with success probability $p$. Let $\tau$ be the minimum of the first jump times of the two processes $X_{z+1}^{(0)}$ and $X_{z+1}^{(1)}$, and define a new process
\be
	X(t) = \begin{cases}
			X_{z+1}^{(1)} & \text{if $\tau > t$,}\\[1ex]
			X_{z+1}^{(0)} & \text{if $\tau \leq t$,}
	       \end{cases}
\ee
describing the well-known Doeblin coupling. Then, let
\ben
	Y_z  = \sum_{k=1}^z X_k^{(1)} + \sum_{k=z+1}^n X_k^{(0)},\qquad Y_{z+1} = Y_z - X_{z+1}^{(0)} + X(t), 			\label{51}
\ee
and one proves that $Y_z$ and $Y_{z+1}$ are Markov processes with generator \eq{48}. Hence, we can write \eq{49} as
\ben
	-\Delta\psi(z) = \int_0^\infty \e^{-t}\IE\bklg{\D h(Y_{z})}\,dt,												\label{52}
\ee
since $\tau$ is exponentially distributed with rate $1$. The bound $\norm{g}\leq 1$ is now immediate from \eq{52}, thus we may assume that $n p q>1$. Note that, from \eq{50} and \eq{51},
\be
	\dist(Y_z) = \Bi(z,p+q\e^{-t})\ast\Bi(n-z,p-p\e^{-t}),
\ee
and hence, from \citet[Lemma~1]{barbour:89},
\besn																												\label{53}
	\V^1\bklr{\dist(Y_z)} & \leq \Var(Y_k)^{-1/2} \\ & \leq \bklr{z(p+q\e^{-t})(q-q\e^{-t})+(n-z)(p-p\e^{-t})(q+p\e^{-t})}^{-1/2} \\
	& \leq \bklr{n p q(1-\e^{-t})}^{-1/2}		.
\ee
Note also that for $\~h := h - 1/2$
\ben
	\babs{\IE\bklg{\D h(Y_{z})}} = \babs{\IE\bklg{\D\~h(Y_{z})}} \leq \V^1\bklr{\dist(Y_z)}/2.													\label{54}
\ee
Thus, applying \eq{54} on \eq{52} and using \eq{53},
\be
	\babs{\D\psi} \leq \int_0^s\e^{-t}\,dt + \frac{1}{2\sqrt{n p q}}\int_s^\infty \frac{\e^{-t}}{\sqrt{1-\e^{-t}}}\,dt.
\ee
Choosing $s = -\ln\bklr{1-(n p q)^{-1}}$ and computing the integrals proves the lemma.
\end{proof}

\subsection{Change of the success probabilities}
\begin{lemma}\label{lemma:a} For every $n\in\IN$, $0<p<1$ and $-(1-p)< t <p$
\ba
	\dtv\bklr{\Bi(n,p-t), \Bi(n,p)} & \leq \abs{t}\bbbklr{
		  \frac{\sqrt{n}}{\sqrt{p q}}
		+ \frac{p-t}{p q} 
		+ \frac{\sqrt{(p-t)(q+t)}}{p q\sqrt{n}}
 		} \\
	\dloc\bklr{\Bi(n,p-t), \Bi(n,p)} & \leq \abs{t}\bbbklr{\frac{1+p-t}{p q} + \frac{\sqrt{(p-t)(q+t)}}{p q\sqrt{n}}} 
\ee
\end{lemma}
\begin{proof} We use Stein's method. If $W\sim\Bi(n,p-t)$, we obtain from \eq{03} and~\eq{04}
\be
	\IE\bklg{(1-p)Wg(W-1)-p(n-W)g(W)} = \IE\bklg{tW\D g(W-1)-tng(W)}
\ee
for every bounded function $g\in F(\IZ)$. The left side is just the Stein operator for $\Bi(n,p)$ hence, taking $g=g_A$ obtained by solving \eq{05} for $\Bi(n,p)$, with the bounds \eq{06} and \eq{47} the $\dtv$-bound follows, noting also that $\IE\abs{W}\leq \abs{\IE W}+\sqrt{\Var W}$. With the remark after \eq{06}, the $\dloc$-bound is proved.
\end{proof}

\subsection{\label{s53}Smoothing properties of independent random variables}

In several parts of this paper, we have the situation that we need to estimate $D^m(U)$, $m=1,2$, for some integer valued random variable $U$, being a sum of some other random variables. If the $U$ is a sum of independent random variables, we can proceed as follows. Assume that $U = \sum_{i=1}^n X_i$, where the $X_i$ are independent. Defining $v_i = \min\klg{\ahalf,1-\ahalf D^1(X_i)}$ and $V = \sum_i v_i$ we obtain from \citet[Proposition~4.6]{barbour:99} the bound
\ben
	D^1(U) \leq \frac{2}{V^{1/2}}																			\label{55}.
\ee
Define further $v^\ast = \max_i  v_i$. Now it is always possible to write $U = U^{(1)} + U^{(2)}$ in such a way that the analogously defined numbers $V^{(1)}$ and $V^{(2)}$ satisfy $V^{(k)} \geq V/2-v^\ast$, $k=1,2$. Using \eq{02} and \eq{55}, we obtain
\ben
	D^2(U) \leq D^1\bklr{U^{(1)}}D^1\bklr{U^{(2)}} \leq \frac{4}{\bklr{V^{(1)} V^{(2)}}^{1/2}} \leq \frac{8}{(V-2v^\ast)_+}.			\label{56}
\ee

\subsection{Smoothing properties of conditional independent random variables}

In most applications, $U$ is a sum of dependent summands and we can not apply \eq{55} and \eq{56} directly. However, assuming that there is a random variable $Z$ on the same probability space as $U$ such that $\dist(U|Z=z)$ can be represented as a sum of independend summands, say $X_i^{(z)}$, $i=1,\dots,n_z$, for each $z$ that $Z$ can attain, we can still apply \eq{55} and \eq{56}, and we obtain
\ban	
	D^1(U) & \leq \IE\bklg{\IE\kle{D^1(U) | Z}} \leq \IE\bbklg{\frac{2}{V_Z^{1/2}}},							\label{57}\\
	D^2(U) & \leq \IE\bklg{\IE\kle{D^1(U) | Z}} \leq \IE\bbklg{\frac{8}{(V_Z-2v_Z^*)_+}},							\label{58}
\ee
where, for each $z$, $V_z$ and $v_z^*$ are the corresponding values as defined in subsection \ref{s53} with respect to the $X_i^{(z)}$.

\section{Appendix}

We now give a generalization of Theorem \ref{th:2}. The proof is omitted, because it runs analogously to the proof of Theorem \ref{th:2}; see also \citet{bkr}.

Suppose that a random variable $W$ satisfies Assumptions G and assume that there are sets $K_i\subset J$, $i\in I$, and square integrable random variables $Z_i$, $Z_{ik}$ and $V_{ik}$, $k\in K_i$ and $i\in I$, as follows:
\bgn
 		\text{$W=W_i + Z_i$, $i\in I$, where $W_i$ is independent of $\xi_i$,}						\label{59}
 	\\
 		Z_i = \sum_{k\in K_i} Z_{ik},																\label{60}
 	\\
 		\text{$W_i = W_{ik} + V_{ik}$, $i\in I$, $k\in K_i$},										\label{61}
 	\\
 		\text{where $W_{ik}$ is independent of the pair $(X_i,Z_{ik})$.}							\notag
\ee
 
\begin{theorem}
With $W$ as above,
\ben
	d_l\bklr{\dist(W),\cBi\bklr{\ceil{4\sigma^2},1/2-t}}\leq \sigma^{-2} \bbklr{\sum_{i\in I}\vartheta_{l,i}+1.75},\qquad l=1,2,	\label{62}
\ee
where
\besn		\label{63}
 	\vartheta_{l,i}  ={} &\frac{1}{2} \IE\bklg{\abs{\xi_i}Z_i^2\,D^l\bklr{\dist(W_i\mid \xi_i,Z_i)}} \\
 			& + \sum_{k\in K_i}\IE\bklg{\babs{\xi_i Z_{ik} V_{ik}}D^l\bklr{\dist(W_{ik}\mid \xi_i,Z_{ik},V_{ik})}}\\
 	& + \sum_{k\in K_i}\babs{\IE\klg{\xi_iZ_{ik}}}\,\IE\bklg{\abs{Z_i +  V_{ik}} D^l\bklr{\dist(W_{ik}\mid Z_i,V_{ik})}}.
\ee
\end{theorem}

\section*{Acknowledgments}

I thank A.~D.~Barbour, D.~Schuhmacher and B.~Nietlispach for many helpful discussions.



\end{document}